\documentclass[12pt,a4paper]{article}

\usepackage[left=2cm,right=2cm,    top=2cm,bottom=2cm,bindingoffset=0cm]{geometry}

\usepackage{subfig}
\usepackage{flushend}
\usepackage{indentfirst}
\usepackage{graphics}
\usepackage{amsmath}
\usepackage{epstopdf}
\usepackage{float}
\usepackage{amssymb}
\usepackage[font=footnotesize,labelfont=bf]{caption}
\usepackage[labelsep=period]{caption}

\newcommand{\be}{\begin{equation}}
\newcommand{\ee}{\end{equation}}

\newtheorem{thm}{Theorem}[section]
\newtheorem{lem}{Lemma}[section]
\newtheorem{nas}{Corrolary}[section]

\newtheorem{ozn}{Definition} [section]

\begin{document}

\title{\textbf{Filtering Problem for Functionals of Stationary Processes with Missing Observations }}

\date{}

\maketitle

\noindent Mathematical Research Press, London, UK\\
Communications in Optimisation Theory, \\
2016, pp.1-18, Article ID 21\\
http://cot.mathres.org

\vspace{20pt}

\author{\textbf{Mykhailo Moklyachuk}$^*$, \textbf{Maria Sidei},\\\\
 {Department of Probability Theory, Statistics and Actuarial
Mathematics, \\
Taras Shevchenko National University of Kyiv, Kyiv 01601, Ukraine}\\
$^{*}$Corresponding Author: Moklyachuk@gmail.com}\\\\\\

\noindent \textbf{Abstract.} \hspace{2pt}
The  problem  of the mean-square optimal linear estimation of the
functional $A\xi=\ \int\limits_{R^s}a(t)\xi(-t)dt,$   which depends on the unknown values of stochastic stationary process $\xi(t)$  from observations of the process $\xi(t)+\eta(t)$ at points $t\in\mathbb{R} ^{-} \backslash S $,  $S=\bigcup\limits_{l=1}^{s}[-M_{l}-N_{l}, \, \ldots, \, -M_{l} ],$ $R^s=[0,\infty) \backslash S^{+},$ $S^{+}=\bigcup\limits_{l=1}^{s}[ M_{l}, \, \ldots, \, M_{l}+N_{l}]$ is considered. Formulas for calculating the mean-square error and the spectral characteristic of the optimal linear estimate of the functional are proposed under the condition of spectral certainty, where spectral densities of the
processes $\xi(t)$ and $\eta(t)$ are exactly known. The minimax (robust)
method of estimation is applied in the case where spectral densities
are not known exactly, but sets of admissible spectral densities are
given. Formulas that determine the least favorable spectral
densities and the minimax spectral characteristics are proposed for some special sets of admissible spectral densities.\\

\noindent \textbf{Keywords:} \hspace{2pt} Stationary process, Mean square error, Minimax-robust estimate, Least favorable spectral density, Minimax spectral characteristic.\\

\noindent \textbf{ Mathematics Subject Classification:} \hspace{2pt}
 Primary: 60G10, 60G25, 60G35, Secondary: 62M20, 93E10, 93E11

\section{Introduction}

Effective methods of solution of the estimation problems (interpolation, extrapolation and filtering) for stationary stochastic sequences and processes were developed by A.~N.~Kolmogorov \cite{Kolmogorov},  N.~Wiener \cite{Wiener}, A.~M.~Yaglom \cite{Yaglom1,Yaglom2}.
An important  contribution to the theory of estimation was made by H. Wold \cite{Wold38,Wold48}, Yu.~A.~Rozanov \cite{Rozanov}, E.~J.~ Hannan \cite{Hannan}.

The basic assumption of most of the developed methods of  estimation of the unobserved values of stochastic processes is that the spectral densities of the considered stochastic processes are exactly known.
However, in practice, these methods can not be applied since the complete information on the spectral densities is impossible in most cases.
In order to solve the estimation problem one have to find  parametric or nonparametric estimates of the unknown spectral densities. Then, under assumption that the selected or estimated densities are the true ones, one of the traditional estimation methods is applied. This procedure can result in significant increasing of the value of error as K.~S.~Vastola and H.~V.~Poor \cite{Vastola} have demonstrated with the help of some examples.
To avoid this effect it is reasonable to search estimates which are optimal for all densities from a certain class of admissible spectral densities. These estimates are called minimax since they minimize the maximum value of the error. This method was  first proposed in  the paper by Ulf ~Grenander \cite{Grenander}  where this approach to extrapolation problem for stationary processes was applied.

Various models of spectral uncertainty and minimax-robust methods of data processing can be found in the survey paper by S.~A.~ Kassam and H.~V.~ Poor \cite{Kassam}. In their papers J. Franke \cite{Franke}, J. Franke and H. V. Poor \cite{Franke_Poor} investigated the minimax extrapolation  and filtering problems for stationary sequences with the help of convex optimization methods. This approach makes it possible to find equations that determine the least favorable spectral densities for some classes of admissible densities.

Papers  by M. Moklyachuk \cite{Moklyachuk:2000} -- \cite{Moklyachuk:2015} contain results of  investigation of the problem of optimal estimation of the functionals which depend on the unknown values of stationary sequences and processes.
M.~ Moklyachuk and A.~ Masyutka  developed the minimax  technique of estimation for vector-valued stationary stochastic processes in papers \cite{Masyutka:2007}--\cite{Masyutka:2012}.
Methods of solution of the  interpolation, extrapolation and filtering problems for periodically correlated stochastic processes were developed by M. Moklyachuk and I. Golichenko  \cite{Golichenko}.
Estimation problems for functionals which depend on the unknown values of stochastic processes with stationary increments were investigated by M. Luz and M. Moklyachuk  \cite{Luz:2014}--\cite{Luz:2015b}.
The interpolation problem  for  stationary sequence with missing values was investigated by M. Moklyachuk and M. Sidei \cite{Sidei1,Sidei2}.

Prediction problem of stationary processes with missing observations  was investigated in papers by P.~Bondon \cite{Bondon1, Bondon2}, Y.~Kasahara, M.~Pourahmadi and A.~Inoue \cite{Kasahara,Pourahmadi},
R. Cheng, A. G. Miamee, M.~Pourahmadi \cite{Cheng}. The interpolation problem for stationary sequences was considered in the paper of H. Salehi \cite{Salehi}.

In this article we deal with the problem of the mean-square optimal linear estimation of the functional $A\xi=\int\limits_ {R^s}a(t)\xi(-t)dt,$  which depends on the unknown values of  a stationary stochastic process $\xi(t)$ from observations of the process $\xi(t)+\eta(t)$ at points  $t\in\mathbb{R}^{-}  \backslash S$,
$S=\bigcup\limits_{l=1}^{s}[ -M_{l}-N_{l},  \ldots,  -M_{l} ],$ $R^s=[0,\infty) \backslash S^{+},$ $S^{+}=\bigcup\limits_{l=1}^{s}[ M_{l}, \, \ldots, \, M_{l}+N_{l}],$ $M_l=\sum\limits_{k=0}^l (N_k+K_k),$   $N_0=0$, $K_0=0$.
The case of spectral certainty as well as the case  of spectral uncertainty are considered. Formulas for calculating the spectral characteristic and the mean-square error of the optimal linear estimate of the functional are derived under the condition that spectral densities of the processes are exactly known. In the case of spectral uncertainty, where the spectral densities are not exactly known while a set of admissible spectral densities is given, the minimax approach is applied. Formulas that determine the least favorable spectral densities and the minimax-robust spectral characteristics of the optimal estimates of the functional are proposed for some specific classes of admissible spectral densities.

\section{Hilbert space projection method of filtering}

Consider a stationary stochastic process $\{\xi(t), t\in \mathbb{R}\}$ with absolutely continuous spectral measure  $F(d\lambda)$ and spectral density $f(\lambda)$. Consider another stationary stochastic process $\{\eta(t), t\in  \mathbb{R}\}$, uncorrelated with the process $\{\xi(t), t\in \mathbb{R}\}$, with absolutely continuous spectral measure  $G(d\lambda)$ and spectral density $g(\lambda)$. Without loss of generality, we  suppose that introduced processes have zero mean values $E\xi(t)=0$, $E\eta(t)=0$.

Assume that the spectral densities  $f(\lambda)$ and $g(\lambda)$ satisfy the minimality condition
\begin{equation}\label{minimal}
\int\limits_{-\infty}^{\infty}\frac{\left|\gamma(\lambda)\right|^2}{f(\lambda)+g(\lambda)}d\lambda<\infty,
\end{equation}
where $\gamma(\lambda)= \int\limits_ {0}^{\infty}\alpha(t)e^{it\lambda}dt$ is a nontrivial function of exponential type.
This condition guarantees that the mean-square errors of estimates of the functionals are nonzero. \cite{Rozanov}.

Stationary processes $\xi(t)$ and $\eta(t)$ admit the spectral decomposition  \cite{Karhunen}
\begin{equation} \label{ksi}
\xi(t)=\int\limits_{-\infty}^{\infty}e^{it\lambda}Z_{\xi}(d\lambda), \hspace{1cm}
\eta(t)=\int\limits_{-\infty}^{\infty}e^{it\lambda}Z_{\eta}(d\lambda),
\end{equation}
where  $Z_{\xi}(d\lambda)$ and  $Z_{\eta}(d\lambda)$ are orthogonal stochastic measures defined on  $[-\pi,\pi)$ that correspond to the spectral measures $F(d\lambda)$ and $G(d\lambda)$, such that the following relations hold true
\[
EZ_{\xi}(\Delta_1)\overline{Z_{\xi}(\Delta_2)}=F(\Delta_1\cap\Delta_2)=\,\frac{1}{2\pi}\int_{\Delta_1\cap\Delta_2}f(\lambda)d\lambda,
\]
\[
EZ_{\eta}(\Delta_1)\overline{Z_{\eta}}(\Delta_2)=G(\Delta_1\cap\Delta_2)=\,\frac{1}{2\pi}\int_{\Delta_1\cap\Delta_2}g(\lambda)d\lambda.
\]

The main purpose  of the article is to find the mean-square optimal linear estimate of the functional
$A\xi= \int\limits_ {R^s}a(t)\xi(-t)dt,$
which depends on the unknown values of the process
$\xi(t)$, based on  the observed values of the process $\xi(t)+\eta(t)$ at points   $t\in\mathbb{R} ^ {-}\backslash S $, where $S=\bigcup\limits_{l=1}^{s}[ -M_{l}-N_{l}, \ldots, -M_{l} ]$, $R^s=[0,\infty) \backslash S^{+},$ $S^{+}=\bigcup\limits_{l=1}^{s}[ M_{l}, \, \ldots, \, M_{l}+N_{l}]$.

Let the function $a(t)$ which determines the functional  $A\xi$ satisfy the conditions
\begin{equation}\label{cond}
\int\limits_{ R^s} \left|a(t)\right|dt<\infty, \quad \int\limits_{R^s} t\left|a(t)\right|^2dt<\infty.
\end{equation}

Due to the spectral decomposition  (\ref{ksi}) of the process $\xi(t)$, the functional $A\xi$ can be represented in the form
\begin{equation*}
A\xi=\int\limits_{-\infty}^{\infty}A(e^{i\lambda})Z_{\xi}(d\lambda), \quad
  A(e^{i\lambda})= \int\limits_ {R^s}a(t)e^{-it\lambda }dt.
 \end{equation*}

Consider the Hilbert space $H=L_2(\Omega,\mathcal{F},P)$ generated by random variables $\xi$ with zero mathematical expectation, $E\xi=0$,  finite variation, $E|\xi|^2<\infty$, and inner product $(\xi,\eta)=E\xi\overline{\eta}$.
Denote by $H^s(\xi+\eta)$ the closed linear subspace generated by elements $\{\xi(t)+\eta(t): t\in \mathbb{R} ^{-} \backslash S\}$ in the Hilbert space $H=L_2(\Omega,\mathcal{F},P)$.
Let  $L_2(f+g)$ be the Hilbert space of complex-valued functions  that are square-integrable  with respect to the measure whose density is  $f(\lambda)+g(\lambda)$, and  $L_2^s(f+g)$ be the subspace of  $L_2(f+g)$ generated by functions $\{e^{it\lambda}, t\in \mathbb{R} ^{-} \backslash S \}.$

Denote by $\hat{A}\xi$ the optimal linear estimate of the functional $A\xi$ from the observations of the process $\xi(t)+\eta(t)$ and by
$\Delta(f,g)=E\left|A\xi-\hat{A}\xi\right|^2$ the mean-square error of the estimate $\hat{A}\xi$.

The mean-square optimal linear estimate $\hat{A}\xi$ of the functional  $A\xi$ is determined by formula
 \begin {equation*}
\hat{A}\xi=\int\limits_{-\infty}^{\infty}h(e^{i\lambda})(Z_{\xi}(d\lambda)+ Z_{\eta}(d\lambda)),
 \end{equation*}
where $h(e^{i\lambda}) \in L_2^s(f+g)$ is the spectral characteristic of the estimate,
and the mean-square error  $\Delta(h;f)$ of the estimate is determined by formula
\begin{equation*}\begin{split}
 \Delta(h;f,g)&=E\left|A\xi-\hat{A}\xi\right|^2=\\
 &=\frac{1}{2\pi}\int\limits_{-\infty}^{\infty}\left|A(e^{i\lambda})-h(e^{i\lambda})\right|^2 f(\lambda)d\lambda
 +\frac{1}{2\pi}\int\limits_{-\infty}^{\infty}\left|h(e^{i\lambda})\right|^2 g(\lambda)d\lambda.
\end{split}\end{equation*}

Since the spectral densities of  stationary processes $\xi(t)$ and $\eta(t)$ are known, in order to find the estimate we can apply the method of orthogonal projections in the Hilbert space proposed by Kolmogorov \cite{Kolmogorov}.
According to this method, the optimal linear estimation of the functional $A\xi$ is a projection of the element $A\xi$ of the space $H$  on  the space $H^s(\xi+\eta)$. The estimate is determined by two conditions:
\begin{equation*} \begin{split}
1)& \hat{A}\xi \in H^s(\xi+\eta), \\
2)& A\xi-\hat{A}\xi \bot  H^s(\xi+\eta).
\end{split} \end{equation*}

Under the second condition the spectral characteristic $h(e^{i\lambda})$ of the optimal linear estimate  $\hat{A}\xi$ for any $t\in \mathbb{R}^{-} \backslash S $ satisfies the relation
\begin{equation*} \begin{split}
&E\left[\left(A\xi-\hat{A}\xi\right)\left(\overline{\xi(t)}+\overline{\eta(t)}\right)\right]=\\
&=\frac{1}{2\pi}\int\limits_{-\infty}^{\infty} \left(A(e^{i\lambda})- h(e^{i\lambda})\right)e^{-it\lambda}f(\lambda)d\lambda-\frac{1}{2\pi}\int\limits_{-\infty}^{\infty}  h(e^{i\lambda})e^{-it\lambda}g(\lambda)d\lambda=0.
\end{split} \end{equation*}

This relation can be written in the following way
\begin{equation}\label{cond2} \begin{split}
\frac{1}{2\pi}\int\limits_{-\infty}^{\infty} \left[A(e^{i\lambda})f(\lambda)- h(e^{i\lambda})(f(\lambda)+g(\lambda))\right]e^{-it\lambda}d\lambda=0, \quad  t\in \mathbb{R}^{-} \backslash S.
\end{split} \end{equation}

Denote the function $C(e^{i\lambda})=A(e^{i\lambda})f(\lambda)- h(e^{i\lambda})(f(\lambda)  +g(\lambda))$, $\lambda \in \mathbb{R}$, and its Fourier transformation
$$\bold{c}(t)=\frac{1}{2\pi}\int\limits_{-\infty}^{\infty} C(e^{i\lambda})e^{-it\lambda}d\lambda,\quad t \in \mathbb{R}.$$

It follows from relation (\ref{cond2}) that the function  $\bold{c}(t)$ is nonzero on the set $T=S \cup [0,\infty)$. Hence,
$$C(e^{i\lambda})=\sum\limits_{l=1}^{s}\int\limits_ {-M_{l}-N_{l}}^{-M_{l}} \bold{c}(t)e^{it\lambda}dt+\int\limits_ {0}^{\infty} \bold{c}(t)e^{it\lambda}dt, $$
and the spectral characteristic of the estimate $\hat{A}\xi$ is of the form
\begin{equation} \label{sphar} \begin{split}
h(e^{i\lambda})=A(e^{i\lambda})\frac{f(\lambda)}{f(\lambda)+g(\lambda)}-\frac{C(e^{i\lambda})}{f(\lambda)+g(\lambda)}.
\end{split} \end{equation}

Under the first condition, $\hat{A}\xi \in H^s(\xi+\eta)$, that determines the estimate of the functional   $A\xi$, for some function $v(t) \in L_2^s(f+g)$ the following relation holds true
$$h(e^{i\lambda})=\frac{1}{2\pi}\int\limits_ {\mathbb{R}^{-}\backslash S}   v( t)e^{it\lambda }d\lambda,$$
therefore, for any  $t\in  T$, we have
\begin{equation}\label{12} \begin{split}
\int\limits_{-\infty}^{\infty}\left(A(e^{i\lambda})\frac{f(\lambda)}{f(\lambda)+g(\lambda)}-\frac{C(e^{i\lambda})}{f(\lambda)+g(\lambda)}\right)e^{-it\lambda}d\lambda=0.
\end{split} \end{equation}

 Let us define the operators in the space $L_2(T)$
 \begin{equation*}\begin{split}
(\bold{B}\bold{x})(t)&=\frac{1}{2\pi}\sum\limits_{l=1}^{s}\int\limits_ {-M_{l}-N_{l}}^{-M_{l}}\bold{x}(u)\int\limits_{-\infty}^{\infty}e^{i\lambda(u-t)}\frac{1}{f(\lambda)+g(\lambda)}d\lambda du+\\
&+\frac{1}{2\pi}\int\limits_ {0}^{\infty}\bold{x}(u)\int\limits_{-\infty}^{\infty}e^{i\lambda(u-t)}\frac{1}{f(\lambda)+g(\lambda)}d\lambda du, \\
(\bold{R}\bold{x})(t)&=\frac{1}{2\pi}\sum\limits_{l=1}^{s}\int\limits_ {-M_{l}-N_{l}}^{-M_{l}}\bold{x}(u)\int\limits_{-\infty}^{\infty}e^{-i\lambda(u+t)}\frac{f(\lambda)}{f(\lambda)+g(\lambda)}d\lambda du+
\\&+\frac{1}{2\pi}\int\limits_ {0}^{\infty}\bold{x}(u)\int\limits_{-\infty}^{\infty}e^{i\lambda(u-t)}\frac{f(\lambda)}{f(\lambda)+g(\lambda)}d\lambda du,\\
(\bold{Q}\bold{x})(t)&=\frac{1}{2\pi}\sum\limits_{l=1}^{s}\int\limits_ {-M_{l}-N_{l}}^{-M_{l}}\bold{x}(u)\int\limits_{-\infty}^{\infty}e^{i\lambda(u-t)}\frac{f(\lambda)g(\lambda)}{f(\lambda)+g(\lambda)}d\lambda du+\\
&+\frac{1}{2\pi}\int\limits_ {0}^{\infty}\bold{x}(u)\int\limits_{-\infty}^{\infty}e^{i\lambda(u-t)}\frac{f(\lambda)g(\lambda)}{f(\lambda)+g(\lambda)}d\lambda du,\\
&\bold{x}(t) \in L_2(T), \quad t \in T.
\end{split}\end{equation*}

The equality (\ref{12}) can be represented in the form
\begin{equation}\label{13} \begin{split}
&\int\limits_{-\infty}^{\infty}\int\limits_{R^s}a(u)e^{i(u-t)}\frac{f(\lambda)}{f(\lambda)+g(\lambda)}dud\lambda\\ &-\left(\int\limits_{-\infty}^{\infty}\left(\sum\limits_{l=1}^{s}\int\limits_ {-M_{l}-N_{l}}^{-M_{l}} \frac{\bold{c}(e^{i(u-t)\lambda})}{f(\lambda)+g(\lambda)}du\right)d\lambda+\int\limits_{-\infty}^{\infty}\int\limits_ {0}^{\infty}\frac{\bold{c}(e^{i(u-t)\lambda})}{f(\lambda)+g(\lambda)}dud\lambda\right)=0.
\end{split} \end{equation}

Denote by  $\bold{a}(t)$  the function such that  $$\bold{a}(t)=0, \, t \in S,\quad \bold{a}(t)=a(t), \, t \in R^s\quad\bold{a}(t)=0, \, t \in S^{+}.$$

Making use of the introduced above denotation, we can represent the equality (\ref{13}) in terms of linear operators in the space $L_2(T)$
\begin{equation}\label{rivn2}
(\bold{R}\bold{a})(t)=(\bold{B}\bold{c})(t), \quad t \in T.
\end{equation}

Assume that the operator $\bold{B}$ is invertible. In this case the function   $\bold{c}(t)$ can be found  by the formula
$$\bold{c}(t)=(\bold{B}^{-1}\bold{R}\bold{a})(t), \quad t \in T.$$

The spectral characteristic $h(e^{i\lambda})$ of the estimate  $\hat{A}\xi$ can be calculated by the formula
\begin{equation}\label{4}
h(e^{i\lambda})=A(e^{i\lambda})\frac{f(\lambda)}{f(\lambda)+g(\lambda)}-
\frac{C(e^{i\lambda}) }{f(\lambda)+g(\lambda)},
 \end{equation}
 \[
C(e^{i\lambda})=\sum\limits_{l=1}^{s}\int\limits_ {-M_{l}-N_{l}}^{-M_{l}} (\bold{B}^{-1}\bold{R}\bold{a})(t)e^{it\lambda}dt+\int\limits_ {0}^{\infty} (\bold{B}^{-1}\bold{R}\bold{a})(t)e^{it\lambda}dt.
\]

The mean-square error of the estimate $\hat{A}\xi$  can be calculated by the formula
\begin{equation} \label{55} \begin{split}
\Delta(h;f,g)=E\left|A\xi-\hat{A}\xi\right|^2&=\frac{1}{2\pi}\int\limits_{-\infty}^{\infty}\frac{\left|A(e^{i\lambda})g(\lambda)+
C(e^{i\lambda}) \right|^2}{(f(\lambda)+g(\lambda))^2}f(\lambda)d\lambda\\
&+\frac{1}{2\pi}\int\limits_{-\infty}^{\infty}\frac{\left|A(e^{i\lambda})f(\lambda)-
C(e^{i\lambda}) \right|^2}{(f(\lambda)+g(\lambda))^2}g(\lambda)d\lambda\\
&=\langle\bold{R}\vec{\bold{a}},\bold{B}^{-1}\bold{R}\vec{\bold{a}}\rangle+\langle\bold{Q}\vec{\bold{a}},\vec{\bold{a}}\rangle,
\end{split} \end{equation}
where $\langle A,C \rangle = \sum\limits_{l=1}^{s}\int\limits_ {-M_{l}-N_{l}}^{-M_{l}}A(t)\overline{C(t)}dt+\int\limits_ {0}^{\infty}A(t)\overline{C(t)}dt$\, is the inner product in the space  $L_2(T)$.

The obtained results can be summarized in the form of theorem.

\begin{thm}\label{t2}
Let  $\xi(t)$ and $\eta(t)$ be uncorrelated stationary processes with spectral densities $f(\lambda)$ and $g(\lambda)$ which satisfy the minimality condition  (\ref{minimal}). The spectral characteristic   $h(e^{i\lambda})$ and the mean-square error  $\Delta(f,g)$ of the optimal linear estimate of the functional  $A\xi$ which depends on the unknown values of the process  $\xi(j)$ based on observations of the process  $\xi(t)+\eta(t),$ $t\in \mathbb{R}^{-}\backslash S$ can be calculated by formulas (\ref{4}), (\ref{55}).
\end{thm}

Let us introduce the notations $N^s = [0,N]\cap R^s$, $R^s=[0,\infty) \backslash S^{+},$ $S^{+}=\bigcup\limits_{l=1}^{s}[ M_{l}, \, \ldots, \, M_{l}+N_{l}].$ Consider the filtering problem for the functional
$A_N\xi= \int\limits_ {N^s}a(t)\xi(-t)dt,$
which depends on the unknown values of the process  $\xi(t)$ based on observations of the process $\xi(t)+\eta(t)$ at time points $t\in\mathbb{R}^{-}\backslash S$.

The optimal linear estimate $\hat{A}_N\xi$  of the functional  $A_N\xi$ is of the form
 \begin {equation*}
\hat{A}_N\xi=\int\limits_{-\infty}^{\infty}h_N(e^{i\lambda})(Z_{\xi}(d\lambda)+ Z_{\eta}(d\lambda)),
 \end{equation*}
where $h_N(e^{i\lambda}) \in L_2^s(f+g)$ is the spectral characteristic of the estimate.

Consider the function $\bold{a}_N(t)$ such that
$$ \bold{a}_N(t)=a(t), \, t \in T\cap [0,N],\quad \bold{a}_N(t)=0, \, t \in T \backslash [0,N].$$

Then the spectral characteristic $h_N(e^{i\lambda})$ and the mean-square error  $\Delta(h_N;f,g)$ of the estimate  $\hat{A}_N\xi$  can be calculated by formulas
\begin{equation}\label{sp_n}
h_N(e^{i\lambda})=A_N(e^{i\lambda})\frac{f(\lambda)}{f(\lambda)+g(\lambda)}-
\frac{C_N(e^{i\lambda}) }{f(\lambda)+g(\lambda)},
 \end{equation}
\[
C_N(e^{i\lambda})=\sum\limits_{l=1}^{s}\int\limits_ {-M_{l}-N_{l}}^{-M_{l}} (\bold{B}^{-1}\bold{R}\bold{a}_N)(t)e^{it\lambda}dt+\int\limits_ {0}^{\infty} (\bold{B}^{-1}\bold{R}\bold{a}_N)(t)e^{it\lambda}dt,
\]
\begin{equation*}
  A_N(e^{i\lambda})= \int\limits_ {N^s}a(t)e^{-it\lambda }dt,
 \end{equation*}
\begin{equation} \label{err_n} \begin{split}
\Delta(h_N;f,g)=E\left|A_N\xi-\hat{A}_N\xi\right|^2=\langle\bold{R}\vec{\bold{a}_N},\bold{B}^{-1}\bold{R}\vec{\bold{a}_N}\rangle+\langle\bold{Q}\vec{\bold{a}_N},\vec{\bold{a}_N}\rangle.
\end{split} \end{equation}

Thus, we obtain the following corollary.

\begin{nas}\label{t2}
Let  $\xi(t)$ and $\eta(t)$ be uncorrelated stationary processes with the spectral densities $f(\lambda)$ and $g(\lambda)$ which satisfy the minimality condition (\ref{minimal}). The spectral characteristic $h_N(e^{i\lambda})$ and the mean-square error $\Delta(h_N;f,g)$ of the optimal linear estimate of the functional  $A_N\xi$ which depends on the unknown values of the process  $\xi(t)$ based on observations of the process  $\xi(t)+\eta(t),$ $t\in \mathbb{R}^{-}\backslash S$ can be calculated by formulas  (\ref{sp_n}), (\ref{err_n}).
\end{nas}

\section{Minimax method of filtering}

 The derived formulas for calculating the value of the mean-square error
 and the spectral characteristic
of the optimal estimate $\hat{A}\xi$ of the functional $A\xi$
can be applied only in the case where we know spectral densities of the processes. However, usually we do not have exact values of the spectral densities of the processes while  we often know a set of admissible spectral densities.
In this case we can apply the minimax-robust approach to estimate the functional $A\xi$.
This method gives us a procedure of finding estimates which minimize the maximum values of the mean-square errors of the estimates for all spectral densities from a given class of admissible spectral
densities simultaneously (see book \cite{Moklyachuk:2008} for more details).

 \begin{ozn}
 For a given class of spectral densities $D=D_f \times D_g$ the spectral densities  $f_0(\lambda) \in D_f$, $g_0(\lambda) \in D_g$ are called the least favorable in the class $D$ for the optimal linear filtering of the functional $A\xi$  if the following relation holds true $$\Delta\left(f_0,g_0\right)=\Delta\left(h\left(f_0,g_0\right);f_0,g_0\right)=\max\limits_{(f,g)\in D_f\times D_g}\Delta\left(h\left(f,g\right);f,g\right).$$
\end{ozn}

\begin{ozn}
For a given class of spectral densities $D=D_f \times D_g$ the spectral characteristic $h^0(e^{i\lambda})$ of the optimal linear filtering  of the functional $A\xi$ is called minimax-robust if there are satisfied conditions
$$h^0(e^{i\lambda})\in H_D= \bigcap\limits_{(f,g)\in D_f\times D_g} L_2^s(f+g),$$
$$\min\limits_{h\in H_D}\max\limits_{(f,g)\in D}\Delta\left(h;f,g\right)=\max\limits_{(f,g)\in D}\Delta\left(h^0;f,g\right).$$
\end{ozn}

From the introduced definitions and formulas derived above we can obtain the following statement.

\begin{lem}
Spectral densities $f_0(\lambda)\in D_f,$ $g_0(\lambda) \in D_g$ satisfying the minimality condition (\ref{minimal}) are the least favorable in the class $D=D_f\times D_g$ for the optimal linear filtering of the functional $A\xi,$ if the Fourier coefficients of the functions $$(f_0(\lambda)+g_0(\lambda))^{-1},\quad f_0(\lambda)(f_0(\lambda)+g_0(\lambda))^{-1}, \quad f_0(\lambda)g_0(\lambda)(f_0(\lambda)+g_0(\lambda))^{-1}$$ form operators $\bold{B}^0, \bold{R}^0, \bold{Q}^0$, which determine a solution to the constrain optimization problem
\begin{equation} \label{extrem} \begin{split}
\max\limits_{(f,g)\in D_f\times D_g}\langle\bold{R}\vec{\bold{a}},\bold{B}^{-1}\bold{R}\vec{\bold{a}}\rangle &+\langle\bold{Q}\vec{\bold{a}},\vec{\bold{a}}\rangle= \\
&\langle\bold{R}^0\vec{\bold{a}},(\bold{B}^0)^{-1}\bold{R}^0\vec{\bold{a}}\rangle+\langle\bold{Q}^0\vec{\bold{a}},\vec{\bold{a}}\rangle.
\end{split}\end{equation}
The minimax spectral characteristic  $h^0=h(f_0,g_0)$ is determined by formula  (\ref{4}) if $h(f_0,g_0) \in H_D.$
\end{lem}

The least favorable spectral densities $f_0(\lambda)$, $g_0(\lambda)$ and the minimax spectral characteristic $h^0=h(f_0,g_0)$ form a saddle point of the function $\Delta \left(h;f,g\right)$ on the set $H_D\times D.$ The saddle point inequalities
$$\Delta\left(h^0;f,g\right)  \leq\Delta\left(h^0;f_0,g_0\right)\leq \Delta\left(h;f_0,g_0\right), \quad \forall h \in H_D, \forall f \in D_f, \forall g \in D_g,$$
hold true if $h^0=h(f_0,g_0)$ and $h(f_0,g_0)\in H_D,$  where $(f_0,g_0)$ is a solution to the constrained optimization problem
\begin{equation} \label{7}
\sup\limits_{(f,g)\in D_f\times D_g}\Delta\left(h(f_0,g_0);f,g\right)=\Delta\left(h(f_0,g_0);f_0,g_0\right),
\end{equation}
\begin{equation*}\begin{split}
\Delta\left(h(f_0,g_0);f,g\right)&=\frac{1}{2\pi}\int\limits_{-\infty}^{\infty}\frac{\left|A(e^{i\lambda})g_0(\lambda)
+C^0(e^{i\lambda})\right|^2}{(f_0(\lambda)+g_0(\lambda))^2}f(\lambda)d\lambda\\
&+\frac{1}{2\pi}\int\limits_{-\infty}^{\infty}\frac{\left|A(e^{i\lambda})f_0(\lambda)-
C^0(e^{i\lambda})\right|^2}{(f_0(\lambda)+g_0(\lambda))^2}g(\lambda)d\lambda,
\end{split}\end{equation*}
\begin{equation*}
C^0(e^{i\lambda})=\sum\limits_{l=1}^{s}\int\limits_ {-M_{l}-N_{l}}^{-M_{l}} ((\bold{B}^0)^{-1}\bold{R}^0\bold{a})(t)e^{it\lambda}dt+\int\limits_ {0}^{\infty} ((\bold{B}^0)^{-1}\bold{R}^0\bold{a})(t)e^{it\lambda}dt, \quad   t \in S.
\end{equation*}

The constrained optimization problem  (\ref{7}) is equivalent to the unconstrained optimization problem \cite{Pshenichnyj}:
\begin{equation} \label{8}
\Delta_D(f,g)=-\Delta(h(f_0,g_0);f,g)+\delta((f,g)\left|D_f\times D_g\right.)\rightarrow \inf,
\end{equation}
 where $\delta((f,g)\left|D_f\times D_g\right.)$ is the indicator function of the set $D=D_f\times D_g$. Solution of the problem (\ref{8}) is characterized by the condition $0 \in \partial\Delta_D(f_0,g_0),$ where $\partial\Delta_D(f_0)$ is the subdifferential of the  functional $\Delta_D(f,g)$ at point $(f_0,g_0)$ \cite{Rockafellar}.

The form of the functional $\Delta(h(f_0,g_0);f,g)$  admits finding derivatives and differentials of the functional in the space $L_1\times L_1$. Therefore the complexity  of the optimization problem (\ref{8}) is determined by the complexity of calculating the subdifferential of the indicator functions  $\delta((f,g)|D_f\times D_g)$  of the sets $D_f\times D_g$ \cite{Ioffe}.

We have the following statement.

\begin{lem}
Let $(f_0,g_0)$ be a solution to the optimization problem (\ref{8}). The spectral densities  $f_0(\lambda)$, $g_0(\lambda)$ are the least favorable in the class $D=D_f\times D_g$ and the spectral characteristic  $h^0=h(f_0,g_0)$ is the minimax of the optimal linear estimate of the functional  $A\xi$ if  $h(f_0,g_0) \in H_D$.
\end{lem}

\section{Least favorable spectral densities in the class $D=D_{\varepsilon_1}^1\times D_{\varepsilon_2}^2$}

Consider the filtering  problem for the functional $A\xi$ in the case where spectral densities of the processes are from to the class of admissible spectral densities $D=D_{\varepsilon_1}^1\times D_{\varepsilon_2}^2$, where
\begin{equation*}
D_{\varepsilon_1}^1=\left\{f(\lambda)\left|\frac{1}{2\pi}\int \limits_{-\infty}^{\infty}\left|f(\lambda)-f_1(\lambda)\right|d\lambda\leq\varepsilon_1\right.\right\}
\end{equation*}
is a  "$\varepsilon$-neighbourhood" in the space $L_1$ of a given bounded spectral density  $f_1(\lambda)$,
\begin{equation*}
D_{\varepsilon_2}^2=\left\{g(\lambda)\left|\frac{1}{2\pi}\int \limits_{-\infty}^{\infty}\left|g(\lambda)-g_1(\lambda)\right|^2d\lambda\leq\varepsilon_2\right.\right\}
\end{equation*}
is a "$\varepsilon$-neighbourhood" in the space $L_2$ of a given bounded spectral density $g_1(\lambda)$.

Suppose that the spectral densities $f_0(\lambda) \in D_{\varepsilon_1}^1$, $g_0(\lambda) \in D_{\varepsilon_2}^2$.
Let the functions determined by the formulas
\begin{equation} \label{hf-filtr}
h_f(f_0,g_0)=\frac{\left|A(e^{i\lambda})g_0(\lambda)+
C^0(e^{i\lambda}) \right|^2}{(f_0(\lambda)+g_0(\lambda))^2},
\end{equation}
\begin{equation} \label{hg-filtr}
h_g(f_0,g_0)=\frac{\left|A(e^{i\lambda})f_0(\lambda)-
C^0(e^{i\lambda})\right|^2}{(f_0(\lambda)+g_0(\lambda))^2},
\end{equation}
be bounded.
Then the functional
$$
\Delta(h(f_0,g_0);f,g)=\frac{1}{2\pi}\int\limits_{-\infty}^{\infty}h_f(f_0,g_0)f(\lambda)d\lambda + \frac{1}{2\pi}\int\limits_{-\infty}^{\infty}h_g(f_0,g_0)g(\lambda)d\lambda.
$$
is continuous and bounded in the space $L_1 \times L_1$. Hence, condition $0 \in \partial\Delta_{D}(f_0,g_0)$, where
$$\partial\Delta_{D_{\varepsilon_1}^1 \times D_{\varepsilon_2}^2}(f_0,g_0)=-\partial\Delta(h(f_0,g_0);f_0,g_0)+\partial\delta((f_0,g_0)\left|D_{\varepsilon_1}^1\times D_{\varepsilon_2}^2\right.),$$
implies that the spectral densities  $f_0(\lambda) \in D_{\varepsilon_1}^1$, $g_0(\lambda) \in D_{\varepsilon_2}^2$ satisfy  equations
 \begin{equation} \label{hf1-filtr}
\left|A(e^{i\lambda})g_0(\lambda)+
C^0(e^{i\lambda})\right|=(f_0(\lambda)+g_0(\lambda))\Psi(\lambda)\alpha_1,
\end{equation}
\begin{equation} \label{hg1-filtr}
\left|A(e^{i\lambda})f_0(\lambda)-
C^0(e^{i\lambda})\right|=(f_0(\lambda)+g_0(\lambda))^2(g_0(\lambda)-g_1(\lambda))\alpha_2,
\end{equation}
where  $\left|\Psi(\lambda)\right|\leq 1$ and $\Psi(\lambda)=sign (f_0(\lambda)-f_1(\lambda))$ if $f_0(\lambda)\neq f_1(\lambda)$, and $\alpha_1 \geq 0$, $\alpha_2 \geq 0$ are constants.

Equations (\ref{hf1-filtr}), (\ref{hg1-filtr}), together with the optimization problem (\ref{extrem}) and normality conditions
\begin{equation} \label{n1-filtr}
\frac{1}{2\pi}\int \limits_{-\infty}^{\infty}\left|f(\lambda)-f_1(\lambda)\right|d\lambda=\varepsilon_1,
\end{equation}
\begin{equation} \label{n2-filtr}
\frac{1}{2\pi}\int \limits_{-\infty}^{\infty}\left|g(\lambda)-g_1(\lambda)\right|^2 d\lambda=\varepsilon_2,
\end{equation}
determine the least favorable spectral densities in the class $D$.

\begin{thm}
Let the spectral densities $f_0(\lambda) \in D_{\varepsilon_1}^1$, $g_0(\lambda) \in D_{\varepsilon_2}^2$ satisfy the minimality condition (\ref{minimal}), and functions determined by formulas  (\ref{hf-filtr}), (\ref{hg-filtr}) be bounded. Spectral densities $f_0(\lambda)$, $g_0(\lambda)$ are the least favorable in the class $ D_{\varepsilon_1}^1 \times D_{\varepsilon_2}^2$ for the optimal linear filtering of the functional $A\xi$ if they satisfy equations  (\ref{hf1-filtr})-- (\ref{n2-filtr}) and determine a solution to the optimization problem (\ref{extrem}). The minimax-robust spectral characteristic of the optimal estimate of the functional $A\xi$ is determined by formula (\ref{4}).
\end{thm}

Consider the problem in the case where spectral densities of the processes are from to the class of admissible spectral densities $D=D_{\varepsilon_1}^2\times D_{\varepsilon_2}^2$, where
$$ D_{\varepsilon_1}^2 = \left\{f(\lambda)\left|\frac{1}{2\pi}\int \limits_{-\infty}^{\infty}\left|f(\lambda)-f_1(\lambda)\right|^2d\lambda\leq\varepsilon_1\right.\right\},$$ $$ D_{\varepsilon_2}^2 = \left\{g(\lambda)\left|\frac{1}{2\pi}\int \limits_{-\infty}^{\infty}\left|g(\lambda)-g_1(\lambda)\right|^2d\lambda\leq\varepsilon_2\right.\right\},$$
where $f_1(\lambda)$, $g_1(\lambda)$ are fixed spectral densities.

\begin{thm}
Suppose that spectral densities $f_0(\lambda) \in D_{\varepsilon_1}^2$, $g_0(\lambda) \in D_{\varepsilon_2}^2$, satisfy the minimality condition (\ref{minimal}) and functions determined by  (\ref{hf-filtr}), (\ref{hg-filtr}) are bounded. Spectral densities $f_0(\lambda)$, $g_0(\lambda)$ are the least favorable in the class $ D_{\varepsilon_1}^2 \times D_{\varepsilon_2}^2$ for the optimal linear filtering of the functional  $A\xi$ if they satisfy   equations
 \begin{equation*}
\left|A(e^{i\lambda})g_0(\lambda)+
C^0(e^{i\lambda})\right|=(f_0(\lambda)+g_0(\lambda))^2(f_0(\lambda)-f_1(\lambda))\alpha_1,
\end{equation*}
\begin{equation*}
\left|A(e^{i\lambda})f_0(\lambda)-
C^0(e^{i\lambda})\right|=(f_0(\lambda)+g_0(\lambda))^2(g_0(\lambda)-g_1(\lambda))\alpha_2,
\end{equation*}
with $(\alpha_1 \geq 0, \, \alpha_2 \geq 0)$,   determine a solution to the optimization problem (\ref{extrem}), and satisfy conditions
\begin{equation*}
\frac{1}{2\pi}\int \limits_{-\infty}^{\infty}\left|f(\lambda)-f_1(\lambda)\right|^2d\lambda=\varepsilon_1,
\end{equation*}
\begin{equation*}
\frac{1}{2\pi}\int \limits_{-\infty}^{\infty}\left|g(\lambda)-g_1(\lambda)\right|^2d\lambda=\varepsilon_2.
\end{equation*}
The function calculated by formula (\ref{4}) is the minimax-robust spectral characteristic of the estimate of the functional $A\xi$.
\end{thm}

\begin{nas}
Assume that the spectral density $g(\lambda)$ is known and the spectral density $f_0(\lambda) \in D_{ \varepsilon_1}^2$. Let the function $f_0(\lambda)+g(\lambda)$ satisfy the minimality condition (\ref{minimal}), and the function $h_f(f_0,g)$ determined by formula (\ref{hf-filtr}) be bounded. The spectral density $f_0(\lambda)$ is the least favorable in the class $D_{ \varepsilon_1}^2$ for the optimal linear filtering of the functional $A\xi$ if it satisfies the relation
\begin{equation*}
\left|A(e^{i\lambda})g(\lambda)+C^0(e^{i\lambda})\right|=(f_0(\lambda)+g(\lambda))^2(f_0(\lambda)-f_1(\lambda))\alpha_1,
\end{equation*}
and the pair $(f_0(\lambda), g(\lambda))$ is a solution of the optimization problem  (\ref{extrem}). The minimax-robust spectral characteristic of the optimal estimate of the functional $A\xi$ is determined by formula (\ref{4}).
\end{nas}

\section{Least favorable spectral densities in the class $D=D_{\varepsilon_1} \times D_{\varepsilon_2}^1$}

Consider the filtering problem  for the functional $A\xi$ in the case where spectral densities of the processes  belong to the class of admissible spectral densities $D_{\varepsilon_1} \times D_{\varepsilon_2}^1$,
$$ D_{\varepsilon_1} = \left\{f(\lambda)\left| f(\lambda)=(1-\varepsilon_1)f_1(\lambda)+\varepsilon_1 w(\lambda),\frac{1}{2\pi}\int\limits_{-\infty}^{\infty} f(\lambda)d\lambda\leq P_1\right.\right\},$$
\begin{equation*}
D_{\varepsilon_2}^1=\left\{g(\lambda)\left|\frac{1}{2\pi}\int \limits_{-\infty}^{\infty}\left|g(\lambda)-g_1(\lambda)\right|d\lambda\leq\varepsilon_2\right.\right\},
\end{equation*}
where spectral densities $f_1(\lambda), g_1(\lambda)$ are fixed, $ w(\lambda)$ is an unknown spectral density. The set $ D_{\varepsilon_1}$ describes the "$\varepsilon$-contamination" \, model of stochastic processes.

Let the spectral densities $f^0(\lambda) \in D_{\varepsilon_1}$, $g^0(\lambda) \in D_{\varepsilon_2}^1$ determine bounded functions  $h_f(f_0,g_0)$, $h_g(f_0,g_0)$ by formulas (\ref{hf-filtr}), (\ref{hg-filtr}). It follows from the condition $0 \in \partial\Delta_{D}(f_0,g_0)$ that the least favorable spectral densities satisfy equations
\begin{equation} \label{hf3-filtr}
\left|A(e^{i\lambda})g_0(\lambda)+
C^0(e^{i\lambda})\right|=(f_0(\lambda)+g_0(\lambda))(\varphi(\lambda)+\alpha_1^{-1}),
\end{equation}
\begin{equation} \label{hg3-filtr}
\left|A(e^{i\lambda})f_0(\lambda)-
C^0(e^{i\lambda})\right|^2=(f_0(\lambda)+g_0(\lambda))\Psi(\lambda)\alpha_2,
\end{equation}
where $\alpha_1, \alpha_2$ are constants,
 $\varphi (\lambda)\leq 0$ and  $\varphi(\lambda)= 0$ if $f_0(\lambda)\geq (1-\varepsilon_1) f_1(\lambda);$ $\left|\Psi(\lambda)\right|\leq 1$ and $\Psi(\lambda)=sign (g_0(\lambda)-g_1(\lambda))$ if $g_0(\lambda)\neq g_1(\lambda)$.

Equations (\ref{hf-filtr}), (\ref{hg-filtr}) together with the extremal condition (\ref{extrem})  and condition
\begin{equation*}
\frac{1}{2\pi}\int \limits_{-\infty}^{\infty}\left|g(\lambda)-g_1(\lambda)\right|d\lambda=\varepsilon_2
\end{equation*}
determine the least favorable spectral densities in the class $D=D_{\varepsilon_1} \times D_{\varepsilon_2}^1$.

The following theorem holds true.

\begin{thm}
Let the spectral densities $f^0(\lambda) \in D_{\varepsilon_1}$, $g^0(\lambda) \in D_{\varepsilon_2}^1$  satisfy the minimality condition (\ref{minimal}), and let functions determined by formulas  (\ref{hf-filtr}), (\ref{hg-filtr}) be bounded. Functions $f_0(\lambda)$, $g_0(\lambda)$ are the least favorable in the class $ D_{\varepsilon_1} \times D_{\varepsilon_2}^1$ for the optimal linear filtering of the functional $A\xi$ if they satisfy equations (\ref{hf3-filtr})-(\ref{hg3-filtr}) and determine a solution to the optimization problem (\ref{extrem}). The function calculated by the formula (\ref{4}) is the minimax-robust spectral characteristic of the optimal estimate of the functional $A\xi$.
\end{thm}

\begin{nas}
Suppose that the spectral density $g(\lambda)$ is known, and the spectral density $f_0(\lambda) \in D_{ \varepsilon_1}$. Let the function $f_0(\lambda)+g(\lambda)$  satisfy the condition (\ref{minimal}), and function $h_f(f_0,g)$ determined by formula (\ref{hf-filtr}) be bounded. The spectral density $f_0(\lambda)$ is the least favorable in the class $D_{ \varepsilon_1}$ for the optimal linear filtering of the functional $A\xi$ if it is of the form
\begin{equation*}
f_0(\lambda)= \max\left\{(1-\varepsilon_1)f_1(\lambda), \alpha_1\left|A(e^{i\lambda})f(\lambda)+
C^0(e^{i\lambda})\right|-g(\lambda)\right\},
\end{equation*}
and the pair  $(f_0(\lambda), g(\lambda))$ is a solution of the optimization problem    (\ref{extrem}). The minimax-robust spectral characteristic of the optimal estimate of the functional $A\xi$ is determined by formula (\ref{4}).
\end{nas}

\begin{nas}
Consider the  known  spectral density  $f(\lambda)$, and the spectral density $g_0(\lambda) \in D_{\varepsilon_2}^1$. Let the function $f(\lambda)+g_0(\lambda)$ satisfy the condition (\ref{minimal}), and function $h_g(f,g_0)$ determined by formula (\ref{hg-filtr}) be bounded. The spectral density  $g_0(\lambda)$ is the least favorable in the class $D_{\varepsilon_2}^1$ for the optimal linear filtering of the functional $A\xi$ if it satisfies the relation
\begin{equation*}
g_0(\lambda)= \max\left\{g_1(\lambda), \alpha_2\left|A (e^{i\lambda})g(\lambda)-
C ^0(e^{i\lambda})\right|-f(\lambda)\right\},
\end{equation*}
and the pair $(f(\lambda), g_0(\lambda))$ determines a solution to the optimization problem (\ref{extrem}). The function calculated by the formula (\ref{4}) is the minimax-robust spectral characteristic of the optimal estimate of the functional $A\xi$.
\end{nas}

\section{Conclusions}

In the article we propose methods of the mean-square optimal linear filtering of the functional which depends on the unknown values of a stationary stochastic process based on observed data of the process with noise and with missing observations. In the case of spectral certainty where spectral densities of the stationary processes are known we derive formulas for calculating the spectral density and the mean-square error of the estimate of the functional.
In the case of spectral uncertainty where spectral densities of the stationary processes are not  exactly known while a set of admissible spectral densities is given, the minimax approach to the filtering problem is applied. Relations that determine the least favorable spectral densities and the minimax-robust spectral characteristics of the optimal estimates of the functional are proposed for some specific classes of admissible spectral densities.

\end{document}